\documentclass[a4paper]{article}
\usepackage{amsmath,amsthm,amsfonts,amssymb}
\title{Families Intersecting on an Interval}
\author{Paul A.~Russell\footnote{Department
of Pure Mathematics and Mathematical Statistics,
Centre for Mathematical Sciences,
Wilberforce Road,
Cambridge CB3 0WB,
England.}
\footnote{{\tt P.A.Russell@dpmms.cam.ac.uk}}}
\date{October 10, 2006}

\newtheorem{theorem}{Theorem}
\newtheorem{lemma}[theorem]{Lemma}

\newcommand{\A}{{\cal A}}
\newcommand{\pow}{{\cal P}}

\newcommand{\sd}{\bigtriangleup}
\newcommand{\B}{{\cal B}}

\newcommand{\modu}{\,\,\,(\!\!\!\!\!\!\mod}
\newcommand{\modv}{\,\,(\!\!\!\!\mod}

\begin{document}
\maketitle

\begin{abstract}
We shall be interested in the following Erd\H os-Ko-Rado-type question.
Fix some set $B\subset[n]$.  How large a family
$\A\subset\pow[n]$ can we find such that the intersection of any two sets
in $\A$ contains a cyclic translate (modulo $n$) of $B$?
Chung, Graham, Frankl and Shearer have proved that, in the case where $B$ is
a block of length $t$, we can do no
better than to take $\A$ to consist of all supersets of $B$.
We give an alternative proof of this result, which is in
a certain sense more `direct'.
\end{abstract}

\begin{section}{Introduction}
Many questions in extremal combinatorics concern the pairwise intersections
of families of subsets of a finite set.  
For example. how large a family $\A\subset[n]^{(r)}
=\{A\subset[n]:|A|=r\}$ can we find with $A\cap B\ne\emptyset$ for all
$A$, $B\in\A$?  This question was answered in the seminal paper of
Erd\H os, Ko and Rado \cite{erdos-ko-rado}:  for $r\le n/2$, we can do
no better than to take $\A$ to be the collection of all $r$-sets containing
some fixed element of $[n]$.  (We note in passing that the question is of
no interest for $r>n/2$ as then the entire family $[n]^{(r)}$ is intersecting.)

Since the publication of \cite{erdos-ko-rado}, the field has rapidly expanded
and is now rich in interesting problems, many of which remain unsolved.
Several such problems arise when we endow the ground-set with some sort
of structure.  
The question we shall be interested in here is the following.  Fix some set 
$B\subset[n]$.  How large a family $\A\subset\pow[n]$ can we find such that
the intersection of any two sets
in $\A$ contains a cyclic translate (modulo $n$) of $B$?
It is conjectured by Chung, Graham, Frankl and Shearer 
\cite{chung-graham-frankl-shearer} that a kernel
system is again best; they are able to establish their conjecture
in the case where $B$ is a block of length $t$:

\begin{theorem}[\cite{chung-graham-frankl-shearer}]\label{block}
Let $n$ and $t$ be positive integers with $t\le n$, and let $\A$ be a family
of subsets of $[n]$ such that whenever we take $A$, $A'\in\A$ then $A\cap A'$
contains some cyclic translate (modulo $n$) of the set $[t]$.  Then 
$|A|\le2^{n-t}$.
\end{theorem}

The aim of this paper is to give an alternative proof of this theorem.  
As we remark
below, it is sufficient to consider instead the problem of finding the
largest possible family of subsets on $[n]$ with any two \emph{agreeing} on
some cyclic translate of the set $[t]$.  The original proof of Theorem
\ref{block} in \cite{chung-graham-frankl-shearer} proceeds in two stages.
First, the authors show that, for $t<n\le 2t$, if $\A$ is a family of subsets
of $[n]$ with any two agreeing on some cyclic translate of $[t]$ \emph{either}
modulo $n$ \emph{or} modulo $n-1$ then $|\A|\le 2^{n-t}$.  They then apply this
result to prove the theorem in general, working in the Abelian group
$(\pow[n],\sd)$  and applying a partitioning argument.  Here, $\sd$ denotes
the symmetric difference operation on the power set $\pow[n]$ of $[n]$, 
i.e.~for $A$, $B\subset[n]$, we define
$$A\sd B=\{x\in A:x\not\in B\}\cup\{x\in B:x\not\in A\}.$$

We again work in the group $(\pow[n],\sd)$.  However, instead of going
via a preliminary result, in our proof we 
show \emph{directly} that $\pow[n]$ can be partitioned into $2^{n-t}$
parts in such a way that no two distinct sets in the same part agree on any
cyclic translate of $[t]$ modulo $n$.
\end{section}

\begin{section}{Algebraic methods}
In this section, we remind the reader of a certain general method for 
bringing algebraic methods
to bear on this sort of problem.  Our problem comes from the general class
of problems of the following form:
\begin{quotation}
\noindent
Suppose $\B$ is some fixed family of subsets of $[n]=\{1,2,\ldots\,,n\}$.
How large can we make a family $\A$ of subsets of $[n]$ subject to the 
condition that for all $A$, $A'\in\A$, there is some $B\in\B$ with 
$B\subset A\cap A'$?
\end{quotation}
We denote by $v(\B)$ the maximal size of a family $\A$ with this property.
Note that the above problems on graphs and arithmetic progressions also fall
into this class.

Unfortunately, the set $\pow[n]$ does not seem to posess any useful algebraic
structure under the intersection operation $\cap$.  However, when endowed
instead with the symmetric difference operation $\sd$, the set $\pow[n]$
becomes an Abelian group.  This leads one to consider a modified version of
the problem, where we insist only that any two sets in $\A$ \emph{agree}
on some set in $\B$:
\begin{quotation}
\noindent
Suppose $\B$ is some fixed family of subsets of $[n]=\{1,2,\ldots\,,n\}$.
How large can we make a family $\A$ of subsets of $[n]$ subject to the 
condition that for all $A$, $A'\in\A$, there is some $B\in\B$ with 
$B\subset[n]-(A\sd A')$?
\end{quotation}
We denote by $\bar{v}(\B)$ the maximal size of a family $\A$ with this 
property.

It is clear that for any family $\B$ we have $v(\B)\le\bar{v}(\B)$.  In
particular, if a kernel family is best for the modified problem then the
same must also be true for the original problem.
Remarkably,
it was proved by Chung, Graham, Frankl and Shearer
that \emph{equality} 
always holds.
This is reassuring, as it means that we know it is always sufficient to attack
the modified problem---a solution to this modified problem will instantly
give a solution to the original problem.

\begin{theorem}[\cite{chung-graham-frankl-shearer}]\label{agree}
Let $\B$ be a family of subsets of $[n]$.  Then $v(\B)=\overline{v}(\B)$.
\end{theorem}

Let us now explain the algebraic idea.
As we have already mentioned, an important advantage of considering our problem
in the modified form above is that, under
the operation $\sd$, $\pow[n]$ forms an Abelian group, and the condition
$B\subset[n]-(A\sd A')$ is equivalent to 
$(A\sd A')\cap B=\emptyset$.  Now,
assume every set in $\B$ has size $t$.  Then we know that $v(\B)\ge 2^{n-t}$
(by considering a kernel system).  Suppose now that we manage to find some
subgroup \hbox{$G\le\pow[n]$} of order $2^t$ 
such that every non-zero set in $G$
intersects every set in $\B$.  Then, given $g\in\pow[n]$ and $h$, $h'\in G$,
we have the set
$(g\sd h)\sd(g\sd h')=h\sd h'$ intersecting every set in $\B$ unless
$h\sd h'=\emptyset$, i.e.~unless $h=h'$.  So any family $\A$ satisfying the
condition that for all $A$, $A'\in\A$ there is some $B\in\B$ with $B\subset
[n]-(A\sd A')$ can contain at most one element from each coset of $G$.  We
may then deduce immediately that $\bar v(\B)\le 2^{n-t}$ and hence that
$v(\B)=2^{n-t}$.

This approach has been used for example by Griggs and Walker
\cite{griggs-walker} to show that if $\B$ consists of all ordinary translates
(rather than cyclic translates) of a fixed set of order $t$ 
then $v(B)=2^{n-t}$,
and by F\"uredi, Griggs, Holzman and Kleitman
\cite{furedi-griggs-holzman-kleitman} to show that if $\B$ consists of all
\emph{cyclic} translates of a fixed set of order 3 then $v(\B)=2^{n-3}$.

In the following section, we apply the method to the case where 
$\B$ consists of all
cyclic translates of a block of length $t$, hence producing a new proof
of Theorem \ref{block}.  The work comes in finding a suitable subgroup $G$,
which in general seems far from obvious.
\end{section}

\begin{section}{Cyclic translates of a block}
We now proceed to our proof of Theorem \ref{block}.  In view of the preceding
section, it is enough to find a subgroup $G$ of $(\pow[n],\sd)$ of order $2^t$
with every non-zero element of $G$ intersecting every block of order $t$.
We shall define the group $G$ by giving a list $g_1$, $g_2$, $\ldots\,$,
$g_t$ of $t$ generators.  For $1\le i,j\le t$, we shall insist that $i\in g_j$
if and only if $i=j$.  This ensures that all of the sums $\sum_{i\in I}g_i$
($I\subset[t]$) are distinct, and hence that $|G|=2^t$.

We begin by considering a number of special cases, beginning with cases
where it is easy to construct the subgroup $G$ and building up to progressively
more complicated cases.  We hope that this will give the reader some feel for
the construction before we come to the (fairly complicated) construction
of $G$ in general.
 
The simplest case of all is where $t|n$.  Then simply take
$$g_i=\{x\in[n]:x\equiv i\modu t)\}.$$
It is clear that each $g_i$ intersects
each block of length $t$.  Moreover, the $g_i$ are pairwise disjoint.  Hence
any non-zero element of $G$ contains some $g_i$, and so intersects every block
of length $t$.

Suppose instead $n\equiv1\modv t)$, say $n=qt+1$.  Then we can take
$$g_i=\{x\in[n]:x\equiv i\modu t)\hbox{ or }x=n\}.$$  By the same reasoning as
above, every non-zero $g\in G$ intersects every block of length $t$ which
is contained entirely within $[n-1]$.  Can some $g$ fail to intersect some
block $B$ of length $t$ with $n\in B$?  If so then $n\not\in g$, and so
$g=\sum_{i\in I}g_i$ for some non-empty $I\subset[t]$ of \emph{even} order.
In particular, $|I|\ge 2$.  Let $i=\min I$ and $j=\max I$.  Then $B$ must
contain at least one of $i$ and $(q-1)t+j$ (as 
$(n+i)-((q-1)t+j)=t+1+i-j\le t$),
and both of these points are in $g$, a contradiction.

More generally, if $n\equiv r\modv t)$ for some $r|t$, say $n=qt+r$,
 then we may take
$$g_i=\{x\in[n]:x\equiv i\modu t)\hbox{ or }(x>n-r\hbox{ and }x\equiv i
\modu r))\}.$$  
The proof that each non-zero $g$ intersects each block of length
$t$ is very similar to the previous case.  The only time when things could
conceivably go wrong is if $g=\sum_{i\in I}g_i$ for some $I\subset[t]$
containing distinct $a$ and $b$ with $a\equiv b\modv r)$.  But then letting
$i$ and $j$ be the least and greatest elements of $I$ congruent to $a$ modulo
$t$, we have $i$, $(q-1)t+j\in g$ and $(n+i)-((q-1)t+j)=t+r+i-j\le t$ and
we are done.

The final special case we consider is where $n\equiv r\modv t)$ for some
$r\nmid t$, say $n=qt+r$, but with $t\equiv r'\modv r)$ for some $r'|r$.  
For $i\le t-r'$ we set
$$g_i=
\{x\in[n]:x\equiv i\modu t)\hbox{ or }(x>qt\hbox{ and }x-qt\equiv i 
\modu r))\}$$
while for $i>t-r'$ we set
$$\{x\in[n]:x\equiv i\modu t)\hbox{ or }(x>qt\hbox{ and }
x-qt\equiv i \modu r'))\}.$$
Again, the $g_i\cap[n-r]$ for $i\in[t]$ are pairwise disjoint, and things
can only go wrong if $g=\sum_{i\in I}g_i$ for some $I\subset[t]$ containing
distinct $i$, $j$ with \hbox{$g_i\cap g_j\ne\emptyset$.}
  There are two ways that
this can happen.  The first is if $I$ contains \hbox{$i\ne j$} with $i\equiv j
\modv r)$,
but we can deal with this case as in the previous paragraph.  The
other possibility is if $I$ contains $i\le t-r'<j$ with 
\hbox{$i\equiv j\modv r')$.}
Let $l$ be the least positive
residue of $i$ modulo $r$, and assume that $i$ is chosen so as to minimize $l$.
We may assume $l>r'$, as otherwise we would have \hbox{$i\equiv j\modv r)$}
 which was
dealt with in our first case.
So $g$ contains each of the points $i$ and
$qt+l-r'$, and $n+i-(qt+l-r')=(i-l)+r+r'\le(t-r-r')+r+r'=t$ so we are done.

We now proceed to define the group $G$ for general $n$ and $t$.  The 
construction can be thought of as an iteration of ideas similar to those used
above.

Let $n$, $t$ be positive integers with $t\le n$.  
We apply Euclid's algorithm to $n$ and $t$, thus obtaining
\begin{eqnarray*}
n&=&q_1t+r_1\\
t&=&q_2r_1+r_2\\
r_1&=&q_3r_2+r_3\\
\vdots&\vdots&\vdots\\
r_{i-2}&=&q_ir_{i-1}+r_i\\
\vdots&\vdots&\vdots\\
r_{k-3}&=&q_{k-1}r_{k-2}+r_{k-1}\\
r_{k-2}&=&q_kr_{k-1},
\end{eqnarray*}
where $t>r_1>r_2>\cdots>r_{k-1}>0$.

Observe that for $k$ odd we have
\begin{eqnarray*}
n&=&q_1t+q_3r_2+q_5r_4+\cdots+q_kr_{k-1}\\
t&=&q_2r_1+q_4r_3+q_6r_5+\cdots+q_{k-1}r_{k-2}+r_{k-1}
\end{eqnarray*}
while for $k$ even we have
\begin{eqnarray*}
n&=&q_1t+q_3r_2+q_5r_4+\cdots+q_{k-1}r_{k-2}+r_{k-1}\\
t&=&q_2r_1+q_4r_3+q_6r_5+\cdots+q_kr_{k-1}.
\end{eqnarray*}
We define the \emph{partial sums} of $n$ by
$$n_m=q_1t+q_3r_2+\cdots+q_{2m-1}r_{2m-2}$$
and of $t$ by
$$t_m=q_2r_1+q_4r_3+\cdots+q_{2m}r_{2m-1},$$
where in each case we allow $m$ to take any value for which the above
expressions make sense.  We interpret
$n_0=t_0=0$ and $n_1=q_1t$.  It will sometimes be convenient to write $r_0=t$
and $r_k=0$.
Observe that we have $n=n_m+r_{2m-1}$ and $t=t_m+r_{2m}$ for each $m$.

Fix $i$ with $1\le i\le t$.
We define $g_i$ in terms of its intersections with the intervals 
$(n_{j-1},n_j]$.
Take $a$ maximal with $t_a<i$.  Now, for $0\le j\le a$, we set
$$g_i^{(j)}=\{x\in(n_j,n_{j+1}]:x-n_j\equiv i-t_{j}\modu r_{2j})\}.$$
So in particular, we have $$g_i^{(0)}=\{1\le x\le n-r_1:x\equiv i\modu t)\}.$$
Define also
$$g_i^{(a+1)}=\left\{\begin{array}{ll}
\{n_{a+1}+[i-t_a]_{r_{2a+1}}\}&\hbox{if }k\ne 2a+1\\
\emptyset&\hbox{if }k=2a+1
\end{array}\right.,$$
where $[y]_z$ denotes the least strictly positive residue of $y$ modulo $z$.
Now, set $g_i=\bigcup_{j=0}^{a+1}g_i^{(j)}$.
We define $G_{n,t}$ to be the subgroup of $\pow[n]$ generated by $g_1$, $g_2$,
$\ldots\,$, $g_t$.  Observe that in the cases $k=1$, 2, 3, this reduces to
our earlier definitions.

\begin{lemma}\label{subgroup}
Let $n$ and $t$ be positive integers with $t\le n$, and define $G_{n,t}$ as
above.  Then
\newcounter{bean}
\begin{list}{(\roman{bean})}{\usecounter{bean}}
\item $|G_{n,t}|=2^t$; and
\item every non-zero element of $G_{n,t}$ intersects every cyclic translate 
modulo $n$ of $[t]$.
\end{list}
\end{lemma}

\begin{proof}
(i) is trivial---observe, for example, that if $1\le s,u\le t$ then $u\in g_s$ 
if and only if $s=u$.

(ii) Let $0\ne g\in G_{n,t}$.  It is enough to show that we can find 
$x_1$,~$x_2$,~$\ldots\,$,~\hbox{$x_m\in g$} with
$x_1<x_2<\cdots<x_m$ satisfying 
$x_{i+1}-x_i\le t$ for $i=1$, 2, $\ldots\,$, $m-1$ and $x_1+n-x_m\le t$.

Suppose first that $g=g_i$ for some $i$ with $1\le i\le t$.  
Then $g$ contains
every $x\in[n]$ with $x\equiv i\modv t)$, so it is enough to show that if
we take $x_1=\min g$ and $x_m=\max g$ then $x_1+n-x_m\le t$.
Now, clearly $x_1=i$.  What is $x_m$?

Take $a$ maximal with $t_a<i$.  If $k\ne 2a+1$ then we must have
\hbox{$x_m=n_{a+1}+[i-t_a]_{r_{2a+1}}$.}
Now, $i-t_a\le t_{a+1}-t_a=q_{2a+2}r_{2a+1}$
and so $(i-t_a)-[i-t_a]_{r_{2a+1}}\le (q_{2a+2}-1)r_{2a+1}$.  Hence
\begin{eqnarray*}
x_1+n-x_m&=&i+n-(n_{a+1}+[i-t_a]_{r_{2a+1}})\\
&=&(i-[i-t_a]_{r_{2a+1}})+(n-n_{a+1})\\
&\le&((q_{2a+2}-1)r_{2a+1}+t_a)+r_{2a+1}\\
&=&t_a+q_{2a+2}r_{2a+1}\\
&\le& t_a+r_{2a}\,=\, t.
\end{eqnarray*}

On the other hand, if $k=2a+1$ then, since $i-t_a\le t-t_a=r_{2a}$, we have 
$x_m=n_{a+1}-r_{2a}+i-t_a=n-r_{2a}+i-t_a$.  Hence
$x_1+n-x_m=i+r_{2a}-i+t_a=t$.

Now, in general, we can write $g=\sum_{i\in I}g_i$ for some non-empty
$I\subset[n]$.  If the $g_i$ ($i\in I$) are pairwise disjoint, then pick
some $i\in I$.  We know that $g_i\subset g$ and that $g_i$ intersects
every block of length $t$.  So $g$ also intersects every block of length $t$.

So we may assume that there exist distinct $i,j\in I$ such that 
$g_i\cap g_j\ne\emptyset$.  
Pick $i$, $j\in I$ with $i<j$ such that $y\in g_i\cap g_j$, 
where $y$ is the least positive 
integer which lies in at least two of the $g_i$ ($i\in I$).  We take
$x_1=i$ and $x_2<\cdots<x_m$ to be the elements of 
$g_j\cap[y-1]$.  As $g_j$ intersects every block of length $t$, it is
enough to check that $n+i-x_m\le t$.

Take $b$ maximal such that $n_b<y$.  Suppose first that $j\le t_b$.  Then,
as $y\in g_j$, we must have $j>t_{b-1}$ and $y=n_b+[j-t_{b-1}]_{r_{2b-1}}$.
Furthermore, $i<j\le t_b$ and $y\in g_i$, so similarly we have
$y=n_b+[i-t_{b-1}]_{r_{2b-1}}$.  Hence $i\equiv j\modv r_{2b-1})$, and, in
particular, $i\le j-r_{2b-1}$.  Now, as $t_{b-1}<j\le t_b$, we know that
$g_j$ contains no elements greater than $n_b$ other than $y$, and that
the elements of $g_j$ in $(n_{b-1},n_b]$ are precisely those 
$x\in(n_{b-1},n_b]$ with $x-n_{b-1}\equiv j-t_{b-1}\modv r_{2b-2})$.  But
$0<j-t_{b-1}\le t-t_{b-1}=r_{2b-2}$ and $r_{2b-2}|n_b-n_{b-1}$.  Hence
$x_m=n_b-r_{2b-2}+j-t_{b-1}$.  So 
\begin{eqnarray*}
i+n-x_m&=&i+n-(n_b-r_{2b-2}+j-t_{b-1})\\
&\le&(j-r_{2b-1})+(n-n_b)+r_{2b-2}-j+t_{b-1}\\
&=&(j-r_{2b-1})+r_{2b-1}+r_{2b-2}-j+(t-r_{2b-2})\\
&=&t,
\end{eqnarray*}
as required.

Now, suppose instead that $j>t_b$.  As $y\in g_j$, we must have 
\hbox{$y-n_b\equiv j-t_b\modv r_{2b})$.}  If we also suppose $i>t_b$ 
then, similarly, we have 
\hbox{$y-n_b\equiv i-t_b\modv r_{2b})$,} and so $i\equiv j\mod r_{2b}$; but
$t\ge i$, $j>t_b=t-r_{2b}$, giving a contradiction as $i\ne j$.  So $i\le t_b$
and $y=n_b+[i-t_{b-1}]_{r_{2b-1}}$.
Now, $j-t_b\le t-t_b=r_{2b}$ so $y\ge j-t_b+n_b$.

If in fact $y=j-t_b+n_b$, then $i-t_{b-1}\equiv j-t_b\modv r_{2b-1})$.  But
$t_{b-1}\equiv t_b\modv r_{2b-1})$, so $i\equiv j\modv r_{2b-1})$ and so
$i\le j-r_{2b-1}$.  Furthermore, $x_m=n_b-r_{2b-2}+j-t_{b-1}$, and so
\begin{eqnarray*}
i+n-x_m&\le&(j-r_{2b-1})+n-(n_b-r_{2b-2}+j-t_{b-1})\\
&=&(j-r_{2b-1})+(n-n_b)+r_{2b-2}-j+t_{b-1}\\
&=&(j-r_{2b-1})+r_{2b-1}+r_{2b-2}-j+(t-r_{2b-2})\\
&=&t,
\end{eqnarray*}
as required.

Otherwise, $y>j-t_b+n_b$. In this case, we have $x_m=y-r_{2b}$.  Now, 
\hbox{$i-t_{b-1}\le t_b-t_{b-1}=q_{2b}r_{2b-1}$.}  So
$(i-t_{b-1})-[i-t_{b-1}]_{r_{2b-1}}\le(q_{2b}-1)r_{2b-1}$, and so
\begin{eqnarray*}
y&\ge&n_b+(i-t_{b-1})-(q_{2b}-1)r_{2b-1}\\
&=&(n_b+r_{2b-1})+i-(t-r_{2b-2})-q_{2b}r_{2b-1}\\
&=&n+i-t+(r_{2b-2}-q_{2b}r_{2b-1})\\
&=&n+i-t+r_{2b}.
\end{eqnarray*}
Hence $i+n-x_m=i+n-(y-r_{2b})\le t$, as required.
\end{proof}

Theorem \ref{block} now follows immediately, as explained earlier.  While it
is interesting to know that Theorem \ref{block} can be proved by this direct
algebraic argument, we cannot at present see any way to generalize this
to deal with cyclic translates of a more general set; the proof seems to
rely heavily on the points of $[t]$ being adjacent.
\end{section}

\end{document}